\documentclass[letterpaper, 10 pt, conference]{ieeeconf}  % Comment this line out if you need a4paper
\newcommand{\fref}[1]{(\ref{eq:#1})}
\usepackage{adjustbox}
\usepackage{array}
\usepackage{booktabs}
\usepackage{multirow}
\usepackage{hhline}

\newcolumntype{R}[2]{%
    >{\adjustbox{angle=#1,lap=\width-(#2)}\bgroup}%
    l%
    <{\egroup}%
}
\IEEEoverridecommandlockouts                              % This command is only needed if 
\usepackage{graphicx}
\usepackage{amsmath}
\usepackage{color}

\title{\LARGE \bf Scheduling of Dynamic Electric Loads Using Energy Storage\\ and Short Term Power Forecasting}

\author{Raymond A. de Callafon, Abdulelah H. Habib and Jan Kleissl
\thanks{R.A. de Callafon, A.H. Habib and J. Kleissl are with the Department of Mechanical and Aerospace Engineering, University of California, San Diego. Email:
       {\tt callafon,ahhabib,jkleissl@ucsd.edu}}}%\\$^*$ Corresponding author}%
\begin{document}

\maketitle
\thispagestyle{empty}
\pagestyle{empty}

%%%%%%%%%%%%%%%%%%%%%%%%%%%%%%%%%%%%%%%%%%%%%%%%%%%%%%%%%%%%%%%%%%%%%%%%%%%%%%%%
\begin{abstract}
In this paper we formulate an optimization approach to schedule electrical loads given a short term prediction of time-varying power production and the ability to store only a limited amount of electrical energy. The proposed approach is unique and versatile as it allows scheduling of electrical loads that each have their own dynamic power demand during on/off switching, while also allowing the specification of minimum on/off times for each loads separately. The optimization approach is formulated as a parallel enumeration of all possible on/off times of the electrical loads using a moving time approach in which only a short term power production forecast is needed, while at the same time taking into account constraints on electrical energy storage and power delivery of a battery system. It is shown that the complexity of the optimization (number of enumerations) is limited by the number of data points in the short term power production forecast and the minimum on/off time of the electrical loads. The limited complexity along with parallel enumeration allows real-time operational scheduling of a large number of loads. The simulation results shown in this paper illustrate that relatively short term power forecast profiles can be used to effectively schedule dynamic loads with various dynamic load profiles. 
\end{abstract}

%%%%%%%%%%%%%%%%%%%%%%%%%%%%%%%%%%%%%%%%%%%%%%%%%%%%%%%%%%%%%%%%%%%%%%%%%%%%%%%%
\section{INTRODUCTION}

One of main challenges in optimal sizing of battery systems for standalone or islanding microgrid applications is the power volatility of distributed renewable energy resources within that microgrid. It is well understood that a combination of battery systems for local energy storage and scheduling of loads for modifying power demands would allow an islanding microgrid to perform reliably, despite power production fluctuations \cite{nrelvar,pvvoltage} and \cite{windsolarchalng}. Even a non-islanding (grid connected) microgrid would benefit from local energy storage and scheduling of loads to allow for power demand following at the point of common coupling to mitigate costly power surge demands and facilitate demand flexibility \cite{An2015}.

Load scheduling and load shedding plays an important role in tracking power production, where loads are turned on/off at optimal times to follow the (predicted) power production as closely as possible \cite{Aponte2006}. Load scheduling and shedding applications can also be used as ancillary services to curtail power surge demands and provide voltage stability. Demonstrations and commercial applications can be seen in EV charging \cite{c2,c4}, building load scheduling \cite{c8}, household appliance scheduling \cite{Loadscheduling}, HVAC system scheduling \cite{c5} and several use cases motivated by environmental and financial incentives \cite{Boon,loadagg,Lang2015932}. An energy storage unit in the form of a battery system can alleviate the effect of power volatility and power demand surges \cite{c1}, but still requires financial decisions with respect to optimal sizing of the battery in combination with load scheduling and load shedding.

Most approaches to optimal load or demand scheduling use some form of Model Predictive Control (MPC) \cite{c12,c13} to compute optimal control or scheduling signals. In MPC typically a constrained optimization problem is solved over a (short term) moving time horizon to compute an optimal control signal for a dynamical system in real time. Countless examples of innovative MPC based approaches for either load scheduling, grid tied storage systems or maintaining voltage stability can be found in \cite{c3}, \cite{c14}, \cite{c15} or \cite{c18}. Although MPC approaches are extremely powerful, models of dynamic systems are typically linear and control signal are allowed to attain any real value during the optimization \cite{c16,c17,c19}. For proper load scheduling and shedding along with a battery, the optimal control may be non-linear and must switch loads with different on/off switching dynamics, while maintaining constraints on battery power and state of charge.

In this paper we define a load scheduling algorithm as the optimal on/off timing of a set of distinct electric loads using an MPC approach in which only a short term power production forecast is needed, while at the same time taking into account constraints on electrical energy storage and power delivery of a battery system via a barrier function. In our problem formulation, each electric load may have a different dynamic power response during on/off switching and a different minimum on/off time. The optimal scheduling is solved by a parallel enumeration of all possible on/off timing combinations that is limited by the number of data points in the short term power production forecast and the minimum on/off time of the electrical loads. The limited complexity along with parallel enumeration of load switching combinations allows real-time operational scheduling of a large number of loads as will be illustrated in the simulation results of this paper.

\section{ELECTRIC LOADS WITH SWITCHING DYNAMICS}\label{sec:SDLM}

\subsection{Switch Signal and Minimum On/Off Time}

To formulate the scheduling algorithm, first the (binary) load switching signal and the transient dynamics describing the (real) power demand of each load is defined here, similar to work presented in \cite{HabibACC2016}. The dynamic properties of a load is characterized by the time-varying power demand during on/off switching and abbreviated to the ``switching dynamics'' of a load. The switching dynamics of each load is important in order to track power production, while maintaining constraints on battery power and state of charge during load switching.

For notational purposes, consider a fixed number of $n$ loads where the time-varying power demand $p_i(t)$, $i=1,2,\ldots,n$ as a function of time $t$ for each load $i$ is determined by a binary switching signal $w_i(t) = \{0,1\}$ used to switch loads ``on'' or ``off''. To maintain generality, the switching dynamics used to describe the power demand $p_i(t)$ during on/off switching may be different for each load, whereas the switching dynamics for a particular load $i$ may also depend on the transition of the binary switching signal $w_i(t) = \{0,1\}$. As a result, the switching dynamics for the time dependent power demand $p_i(t)$ for each load $i$ is different when the binary switching signal $w_i(t) = \{0,1\}$ transitions from 0 to 1 (rising edge, turning load $i$ on) or from 1 to 0 (falling edge, turning load $i$ off). 

Each load $i$ is also assumed to have a known minimum duration $T_i^{\mbox{\tiny on}}>0$ for the ``on'' time of the load when $w_i(t)=1$ and a minimum duration $T_i^{\mbox{\tiny off}}>0$ for the ``off'' time of the load when $w_i(t)=0$. The minimum on/off times $T_i^{\mbox{\tiny off}}$ and $T_i^{\mbox{\tiny on}}$ avoid undesirable chattering of the switch signal $w_i(t)$ during load scheduling and directly reduce the number of possible switching combinations.

For the computation of the time varying binary switch signal $w_i(t)$, a Model Predictive Control (MPC) approach is used over a moving time optimization horizon $T< \infty$. The finite optimization horizon $T$ is required to satisfy 
\begin{equation}
\max_i (T_i^{\mbox{\tiny on}}, T_i^{\mbox{\tiny off}}) < T <\infty
\label{eq:Tcons}
\end{equation}
to ensure the effect of the load switching signal $w_i(\tau)$ at $\tau=t$ for loads with minimum on/off times $T_i^{\mbox{\tiny on}}$ and $T_i^{\mbox{\tiny off}}$ can be predicted into the future over the time interval $T$ when $t \leq \tau \leq t+T$. Finally, it is assumed that all loads are initially switched ``off'', e.g., $w_i(0)=0$ for $i=1,2,\ldots,n$.

\subsection{Finite Set of Admissible Switching Signals}

With the minimum on/off duration times $T_i^{\mbox{\tiny on}},T_i^{\mbox{\tiny off}}$ and the finite optimization horizon $T$ for load switching, on/off switching of a load at time $t$ can now be formalized. The load switching signal $w_i(t),~ i=1,2,\ldots,n$ will be a Zero Order Hold (ZOH) binary signal, where $w_i(t) \in  \{0,1\}$ is allowed to change its binary state only over a finite number $N$ of switching time opportunities within the prediction horizon of length $T$. For binary load switching, an MPC optimization problem quickly becomes intractable due to a combinatorial problem where the number of switching combinations grows exponentially in the number of switching opportunities $N$ and the number $n$ of loads. However, it was shown in the \cite{HabibACC2016} that the requirement of the minimum on/off duration times $T_i^{\mbox{\tiny on}},T_i^{\mbox{\tiny off}}$ over the finite time prediction horizon $T$ constrained by \fref{Tcons} significantly reduces the number of switching combinations and alleviates the combinatorial problem. In fact, it is shown in \cite{HabibACC2016} that the number of possible combinations is much smaller than a trivial exponential growth of $(2^n)^{N-1}$.

With the ZOH approximation, the switching signal $w_i(t) \in  \{0,1\}$ is kept constant in between the finite number $N$ of switching time opportunities within the prediction horizon of length $T$. As a result, the admissible on/off transition signal $w_i(t) = \{w_i^{\mbox{\tiny on}}(t),w_i^{\mbox{\tiny off}}(t)\}$ of a load at time $t=\tau_i$ can now be formalized by the ZOH switching signal
\begin{equation}
w_i^{\mbox{\tiny on}}(t) = 
\left \{ 
\begin{array}{rclcl} 
0 &\mbox{for}& t < \tau_i & \mbox{and}& \tau_i \geq T_{i,last}^{\mbox{\tiny off}} + T_i^{\mbox{\tiny off}} \\ 
1 &\mbox{for}&  t \geq \tau_i & \mbox{and} & \tau_i \leq T - T_i^{\mbox{\tiny on}}
\end{array} \right .
\label{eq:wion}
\end{equation}
where $T_{i,last}^{\mbox{\tiny off}}$ denotes the most recent (last) time stamp at which the load $i$ was switched ``off'', and
\begin{equation}
w_i^{\mbox{\tiny off}}(t) = 
\left \{ 
\begin{array}{rclcl} 
1 &\mbox{for}& t < \tau_i & \mbox{and}& \tau_i \geq T_{i,last}^{\mbox{\tiny on}} + T_i^{\mbox{\tiny on}} \\ 
0 &\mbox{for}&  t \geq \tau_i
\end{array} \right .
\label{eq:wioff}
\end{equation}
where $T_{i,last}^{\mbox{\tiny on}}$ denotes the most recent (last) time stamp at which the load $i$ was switched ``on''.

\subsection{Dynamic Load Models}

It is clear that the switching time(s) $\tau_i$ for the signal $w_i(t) = \{0,1\}$ depends on the time-varying power demand $p_i(t)$ that is different for each load. For the computational results presented in this paper, continuous-time (linear) dynamic models will be used to model the power switching dynamics of a load. It should be pointed out that the computational analysis is not limited to the use of a linear dynamic model, as long as the dynamic model allow the numerical computation of a time-varying power demand $p_i(t)$ as a function of the switching signal $w_i(t)$ over a finite time interval $T$.

To allow different dynamics for the time dependent power demands $p_i(t)$ when the binary switching signal $w_i(t) = \{0,1\}$ transitions from 0 to 1 ("on") or from 1 to 0 ("off"), different dynamics is used for each of the load models. This allows power demands $p_i(t)$ to be modeled at different rates when switching loads. Using the Laplace transform ${\cal L}\{ \cdot \}$ and referring back to the admissible on/off transition signals $w_i^{\mbox{\tiny on}}(t)$ and $w_i^{\mbox{\tiny off}}(t)$ respectively in \fref{wion} and \fref{wioff}, the switched linear order continuous-time dynamic models for the loads are assumed to be of the form
\begin{equation}
p_i(s) = G_{i}^{\mbox{\tiny on}}(s) x_i w_i(s)~\mbox{and}~ w_i(s) = {\cal L} \{ w_i^{\mbox{\tiny on}}(t)\}
\label{eq:odeon}
\end{equation}
and
\begin{equation}
p_i(s) = G_{i}^{\mbox{\tiny off}}(s) x_i w_i(s)~\mbox{and}~ w_i(s) = {\cal L} \{ w_i^{\mbox{\tiny off}}(t)\}
\label{eq:odeoff}
\end{equation}
where $G_{i}^{\mbox{\tiny on}}(s)$ and $G_{i}^{\mbox{\tiny on}}(s)$ represent the dynamics of the power demand for turning the load $i$ "on" or "off". Both models satisfy $G_{i}^{\mbox{\tiny on}}(0)=1$ and $G_{i}^{\mbox{\tiny on}}(0)=1$ and a steady-state load demand parameter $x_i$ is used to model the relative size of the load, but different dynamics is used to model respectively the on/off dynamic switching of the load \cite{HabibACC2016}. 

\subsection{Discretization of Switching Dynamics}

In order to be able to compute the time dependent power demand $p_i(t)$ for each load $i$, the time response of the switched dynamic models given in \fref{odeon} and \fref{odeoff} needs to be computed over the prediction horizon $T$. The switching signal $w_i(t)$ and the power demand $p_i(t)$ for each load $i$ is time discretized at $t_k = k \Delta_t$ where $\Delta_t$ is the sampling time $k=0,1,\ldots$ is an integer index.

The discretized load switching signal $w_i(t_k),~ i=1,2,\ldots,n$ is a Zero Order Hold (ZOH) binary signal. Since $w_i(t_k) \in  \{0,1\}$ is allowed to change its binary state only over a finite number $N$ of switching times $\tau_i$ within the prediction horizon of length $T$, we assume that 
both the switching times
\begin{equation}
\tau_i = N_i \Delta_t
\label{eq:Ni}
\end{equation}
and the minimum on/off duration times
\begin{equation}
\begin{array}{rcl}
 T_i^{\mbox{\tiny on}} &=& N_i^{\mbox{\tiny on}} \Delta_t \\
 T_i^{\mbox{\tiny off}} &=& N_i^{\mbox{\tiny off}} \Delta_t  
\end{array}
\label{eq:Nionoff}
\end{equation}
are all multiple of the sampling time $\Delta_t$.

As the switching signal $w_i(t_k)$ is {\em also\/} held constant between subsequent time samples $t_k$ and $t_{k+1}$, the computation of the time discretized power demand $p_i(t_k)$ for each load can be achieved using a Zero Order Hold (ZOH) discrete-time equivalent of the continuous-time models given earlier in \fref{odeon} and \fref{odeoff}. Using the z-transform ${\cal Z}\{ \cdot \}$, the ZOH discrete-time equivalent dynamic models are given by
\[
p_i(z) = G_{i}^{\mbox{\tiny on}}(z) x_i w_i(z)~\mbox{and}~ w_i(z) = {\cal Z} \{ w_i^{\mbox{\tiny on}}(t_k)\}
\]
for ``on'' switching of the load and 
\[
p_i(z) = G_{i}^{\mbox{\tiny off}}(z) x_i w_i(z)~\mbox{and}~ w_i(z) = {\cal Z} \{ w_i^{\mbox{\tiny off}}(t_k)\}
\]
for ``off'' switching of the load, where $G_{i}^{\mbox{\tiny on}}(z)$ and $G_{i}^{\mbox{\tiny off}}(z)$ are the ZOH discrete-time equivalents of $G_{i}^{\mbox{\tiny on}}(s)$ and $G_{i}^{\mbox{\tiny off}}(s)$ using a sampling time $\Delta_t$. As a result, the power demand dynamics of each load is fully determined by $G_{i}^{\mbox{\tiny on}}(z)$, $G_{i}^{\mbox{\tiny off}}(z)$, static load demand $x_i$ and the chosen sampling time $\Delta_t$. 

\section{LOAD SCHEDULING ALGORITHM}\label{sec:DLS}

\subsection{Battery System Constraints}

The discrete-time switching signals $w_i(t_k)$ for the loads $i=1,2,\ldots,n$ lead to a total (real) power demand
\[
p(t_k) = \sum_{i=1}^n p_i(t_k)
\]
for the $n$ schedulable loads. To formulate a dynamic load scheduling algorithm, first the power tracking error 
\begin{equation}
e(t_k) = P(t_k) -\sum_{i=1}^n p_i(t_k)
\label{eq:powererror}
\end{equation}
is defined as the error between the anticipated or predicted power generation $P(t_k)$ and total real power demand $p(t_k)$ is considered. The time-dependent variable $P(t_k)$ may refer to any result obtained from net power generation and prediction, e.g. solar power prediction. Since power generation prediction is not the main objective or contribution of this paper, $P(t_k)$ is left as general as possible here. Due to the limited number $n$ of loads and possible errors in power prediction $P(t_k)$, any (predicted) power tracking error $e(t_k)$ will be absorbed/delivered by an energy storage system to match net power flow.

Using a battery system for energy storage and assuming power tracking errors $e(t_k)$ can be absorbed or delivered by the battery, the signal $e(t_k)$ is subjected to several constraints imposed by the battery system. For a typical battery system, the first constraint involves the maximum power delivery and absorption capability
\begin{equation}
P \cdot |e(t_k)|< 1
\label{eq:batpowercons}
\end{equation}
of the battery normalized to Power Units (PU) by $P$. The second constraint involves the maximum energy storage capability
\begin{equation}
0.1 \leq S \cdot \sum_{m=0}^k e(t_k) < 0.9 
\label{eq:batenergycons}
\end{equation}
of the battery normalized to State of Charge (SOC) Units by $S$. The lower bound of 10\% and upper of 90\% is chosen as a safeguard to protect the battery against under- and overcharging, but can be chosen closer to 0 or 1 if so desired.

\subsection{Admissible Discrete-Time Switching Combinations}

With the imposed time discretization given in \fref{Ni}, \fref{Nionoff} and a finite number $N$ of switching times $\tau_i$ within the prediction horizon of length $T$, the admissible on/off transition signal in \fref{wion} reduces to
\begin{equation}
w_i^{\mbox{\tiny on}}(t_k) = 
\left \{ 
\begin{array}{rclcl} 
0 &\mbox{for}& k < N_i & \mbox{and}& N_i \geq N_{i,last}^{\mbox{\tiny off}} + N_i^{\mbox{\tiny off}} \\ 
1 &\mbox{for}&  k \geq N_i & \mbox{and} & N_i \leq N - N_i^{\mbox{\tiny on}}
\end{array} \right .
\label{eq:wiondiscrete}
\end{equation}
where $N_{i,last}^{\mbox{\tiny off}}$ now denotes the most recent discrete-time index at which the load $i$ was switched ``off''. Similarly, \fref{wioff} reduces to
\begin{equation}
w_i^{\mbox{\tiny off}}(t_k) = 
\left \{ 
\begin{array}{rclcl} 
1 &\mbox{for}& k < N_i & \mbox{and}& N_i \geq N_{i,last}^{\mbox{\tiny on}} + N_i^{\mbox{\tiny on}} \\ 
0 &\mbox{for}&  k \geq N_i
\end{array} \right .
\label{eq:wioffdiscrete}
\end{equation}
where $N_{i,last}^{\mbox{\tiny on}}$ denotes the most recent discrete-time index at which the load $i$ was switched ``on''. Collectively, the signals $w_i^{\mbox{\tiny on}}(t_k)$ in \fref{wiondiscrete} and $w_i^{\mbox{\tiny off}}(t_k)$ \fref{wioffdiscrete} define a set ${\cal W}$ of binary values for admissible discrete-time switching signals defined by
\begin{equation}
{\cal W} = \left \{ 
\begin{array}{c}
w_i(t_k) \in \{w_i^{\mbox{\tiny on}}(t_k),w_i^{\mbox{\tiny off}}(t_k)\} ,\\
 i=1,2,\ldots,n,~ k=1,2,\ldots,N\\
\mbox{where}~
\begin{array}{c}
w_i^{\mbox{\tiny on}}(t_k) \in \{0,1\}~\mbox{given in}~\mbox{\fref{wiondiscrete}}\\
w_i^{\mbox{\tiny off}}(t_k) \in \{0,1\}~\mbox{given in}~\mbox{\fref{wioffdiscrete}}
\end{array}
\end{array}
\right \}
\label{eq:Wset}
\end{equation}

It is worthwhile to note that the number of binary elements in the set ${\cal W}$ in \fref{Wset} is always (much) smaller than the trivial exponential number of $(2^n)^{N-1}$ \cite{HabibACC2016}. This due to required minimum number of on/off samples $N_i^{\mbox{\tiny on}},N_i^{\mbox{\tiny off}}$ for the loads given in \fref{wiondiscrete} and \fref{wioffdiscrete}.

\subsection{Moving Horizon Optimization}

Following the power tracking error $e(t_k)$ defined in \fref{powererror}, the dynamic load scheduling optimization problem is formulated as a moving horizon optimization problem
\begin{equation}
\begin{array}{c}w_i(t_m) \\ i=1,2,\ldots,n \\ 
m=k,\ldots,k+N-1\end{array} = \mbox{arg} \min_{w_i(t_m) \in {\cal W}}  f(e(t_m)),
\label{eq:MPC}
\end{equation}
where $f(e(t_m))\geq0$, $i=1,2,\ldots,n$ refers to the $n$ loads and $m=k,\ldots,k+N-1$ refers to the $N$ switching time combinations within a prediction horizon of length $T$ over the admissible set ${\cal W}$ defined in \fref{Wset}. Similar to the ideas in Model Predictive Control (MPC), the $N \times n$ dimensional optimal switching signal $w_i(t_m)$ is computed over the prediction horizon $m=k,\ldots,k+N-1$. Once the optimal switching signal $w_i(t_m) \in {\cal W}$, $m=k,\ldots,k+N-1$ is computed, the optimal signal is applied to the loads {\em only\/} at the time instant $t_k$, after which the time index $k$ is incremented and the optimization in \fref{MPC} is recomputed over the moving time horizon. 

It should be noted that the admissible set ${\cal W}$ defined in \fref{Wset} has a finite and countable number of binary combinations for the switching signal that is (much) smaller than the trivial exponential number of $(2^n)^{N-1}$ \cite{HabibACC2016}. Therefore, the $N \times n$ dimensional optimal switching signal $w_i(t_m) \in {\cal W}$ is computed simply by a finite number of evaluation of the criterion function $f(e(t_l))>0$. Furthermore, evaluation of the power tracking error $e(t_m)$ in \fref{powererror} for each possible switching combination $w_i(t_m) \in {\cal W}$ can be done with a full parallel computation, as different $w_i(t_m)$ for $i=1,2,\ldots,n$ and $m=k,2,\ldots,k+N-1$ are independent of each other. Instead of formulating a gradient based optimization or Mixed Integer Linear Programming (MILP) for loads with different switching dynamics, this approach allows extremely fast (parallel) numerical evaluation of the power tracking error $e(t_m)$ in \fref{powererror} and the battery constraints \fref{batpowercons} and \fref{batenergycons} for the finite number of switching signal combinations $w_i(t_m)$ within the set ${\cal W}$ for real-time operation dynamic scheduling of loads.

%It is clear that evaluation of \fref{MPC} does require a predicted power $P(t_m)$ over the prediction horizon $m=k,\ldots,k+N-1$ of length $T$, but the length of the prediction horizon can be kept relatively short. The power prediction horizon $T$ is only subjected to the constraints given earlier in \fref{Tcons} to ensure the effect of the load switching signal $w_i(t_k)$ for loads with a minimum on time $T_i^{\mbox{\tiny on}}$ and off time $T_i^{\mbox{\tiny off}})$ can be simulated into the future over the prediction horizon with length $T$.

To incorporate the battery constraints given earlier in \fref{batpowercons} and \fref{batenergycons} and to be able to track the predicted power $P(t_k)$, the  optimization function $f(e(t_m))>0$ is defined as the sum of a least squares criterion $\| e(t_m) \|_2$ and (smooth) boundary functions $B_j(e(e(t_m))$. In particular, the optimization function $f(e(t_m))$ is defined as
\begin{equation}
f(e(t_m)) = \| e(t_m) \|_2  + \sum_{j=1}^4 B_j(e(t_m))
\label{eq:optcrit}
\end{equation}
where
\[
\| e(t_m) \|_2 = \sum_{m=k+1}^{k+N-1} tr\{e(t_m) e(t_m)^T
\]
and the barrier functions $B_j(e(t_m))$ are defined as follows

\begin{itemize}

\item $\displaystyle B_1(e(t_m)) = C_1 \cdot (P \max_{m} |e(t_m)| -1)$\\ if $\displaystyle \max_{m} P|e(t_m)| \geq 1$, else $B_1(e(t_m))=0$.

\item $\displaystyle B_2(e(t_m)) = - C_2 \cdot \Delta_t \sum_{m} e(t_m)$\\ if $\displaystyle \sum_{m} e(t_m) \leq 0$, else $B_2(e(t_m))=0$.

\item $\displaystyle B_3(e(t_m)) = C_3 \cdot (S \Delta_t \sum_{m} e(t_m) -0.9)$\\ if $\displaystyle \sum_{m} e(t_m) \geq 0.9$, else $B_3(e(t_m))=0$.

\item $\displaystyle B_4(e(t_m)) = - C_4 \cdot (S \Delta_t \sum_{m} e(t_m) - 0.1)$\\ if $\displaystyle \sum_{m} e(t_m) \leq 0.1$, else $B_4(e(t_m))=0$.

\end{itemize}

The elaborate definition of $f(e(t_m))$ in \fref{optcrit} ensures that $f(e(t_m))\geq0$ and the constraints \fref{batpowercons} and \fref{batenergycons} are taken into account with a linear weighting scaled by the constants $C_j$. Although the barrier functions $B_j(e(t_m))$ are not ``true'' barrier functions that approach $\infty$ at the constraint, the additional linear weighting ensures that solutions are found that are forced away from the constraints. Increasing the value of $C_j$ will make this enforcement stronger and typically $C_2 >> 1$ to enforce that the normalized SOC
\[
S \Delta_t \sum_{m} e(t_m) 
\]
always remains positive. In the application example used in this paper, the coefficients $C_j$ were set to $C_1 = C_3 = C_4 = 10$ whereas $C_2 = 1000$. The reason why not a true barrier function is used is that the optimization of the function $f(e(t_m))\geq0$ in \fref{MPC} still allows for a solution in case the constraints are (temporarily) violated instead of giving no possible solution for load scheduling. Temporarily violation of constraints can be used during battery storage design to indicate that a larger battery is required, while in operation it may be used to allow for a (temporary) solution for load scheduling instead of providing simply ``no'' solution due unanticipated constrain violation.

\section{Application Example}\label{sec:AEx}

\subsection{Simulated Load Switching Dynamics}

To illustrate the results of the scheduling algorithm for loads with distinct load switching dynamics, three loads are selected with different switching dynamics for on/off switching. Load sizes were elected using optimal static load size selection \cite{HabibACC2016}. Furthermore, load no. 2 exhibits a resonant power dynamics behavior requiring a temporary power surge to power on the load. The dynamic characteristics of the loads with their minimum on/off time used in this simulation study are summarized in Table~\ref{table1}.
\begin{table}[h]
\centering
\caption{Loads characteristics: relative size in PU, poles of denominator dynamics, and minimum on/off time in seconds.}
\label{table1}
\begin{tabular}{|c|c|c|c|c|c|}
\hline
\begin{tabular}[c]{@{}c@{}}Loads \\ char.\end{tabular} & Size (\%) & Poles$^{\mbox{\tiny on}}_i$             & Poles$^{\mbox{\tiny off}}_i$ & $T_i^{\mbox{\tiny on}}$ & $T_i^{\mbox{\tiny off}}$ \\ \hline
$x_1$                                                         & 60.00        & -0.01                & -0.04     & 180        & 180         \\ \hline
$x_2$                                                         & 25.86       & - 0.05 $\pm$ j0.06 & -0.05     & 240        & 240         \\ \hline
$x_3$                                                         & 12.22        & -0.02              & -0.02   & 300        & 300         \\ \hline
\end{tabular}
\end{table}

\begin{figure}[ht]
\centering
\includegraphics[width=.85\columnwidth]{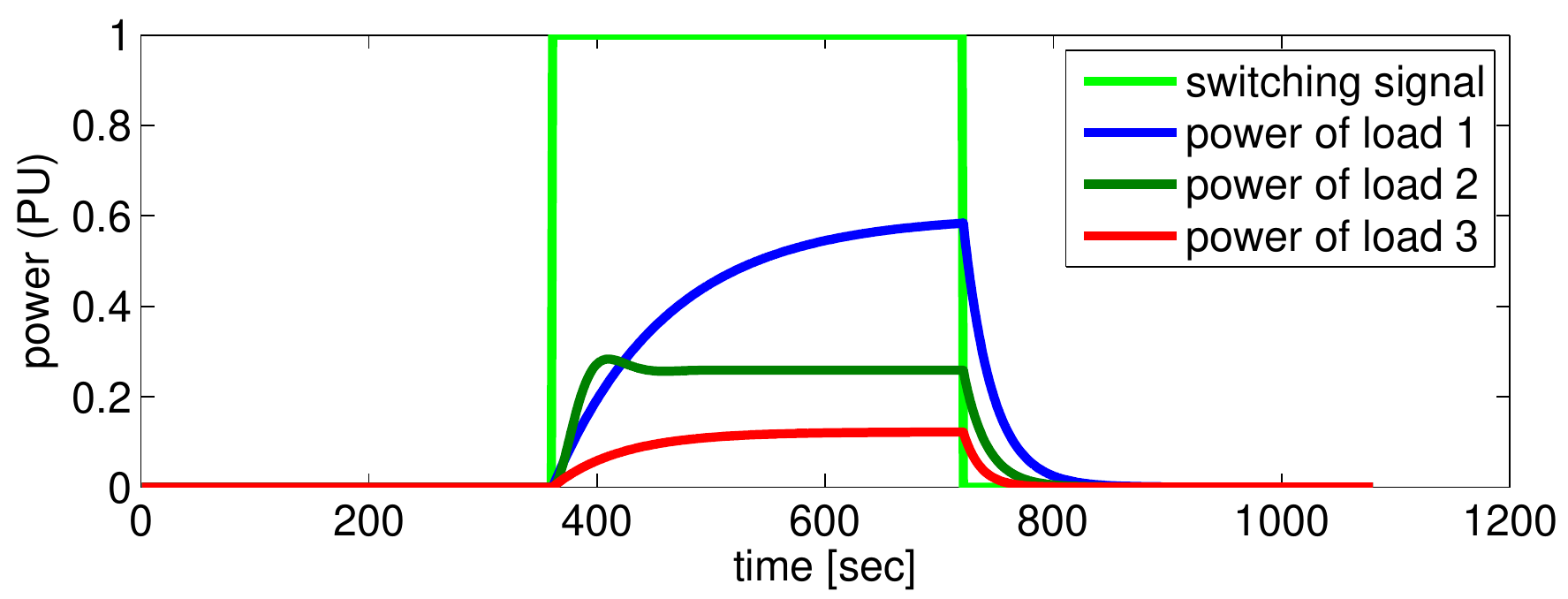}
\caption{Dynamics of power demand $p_i(t_k)$ (colored lines) of the three loads $i=1,2,3$ defined in Table~\ref{table1} as a function of the same binary switch signal $w_i(t_k)$. Note: time scale is in seconds.}
\label{dynamicsload}
\end{figure}

To illustrate the variability in the dynamics of the loads summarized in Table~\ref{table1}, the dynamic response of switching dynamics of the three loads in our case study are depicted in Figure~\ref{dynamicsload}. Although the same switching signal $w_i(t_k)$ is used, the loads exhibit different power demand transitions $p_i(t_k)$. Load 2 shows the typical behavior of a second order dynamics model with an initial larger peak load, typically seen in AC motors used in HVAC systems. 

\subsection{Dynamic Load Switching Results at 50\% SOC}

For the simulation results in this paper, the power curve to be tracked is a power production curve produced by a solar power unit with an irregular bell shaped curve due to solar variability during a 4 hour (240 minute) period. It is worth mentioning that a single optimization for $n=3$ loads over a prediction horizon of $N=6$ switching time opportunities (every minute) with the load dynamics summarized in Table~\ref{table1} (discretized at every second) takes less then 0.3 second to compute in Matlab on a standard 4 core CPU system. The optimization results displayed in each of the figures that follow was therefore computed in less than 70 seconds over the 240 switching opportunities along the 4 hour power curve. As loads are scheduled to switch each minute, the 0.3 second optimization time clearly poses no problem for real-time operation. 

\begin{figure}[ht]
\centering
\includegraphics[width=.85\columnwidth]{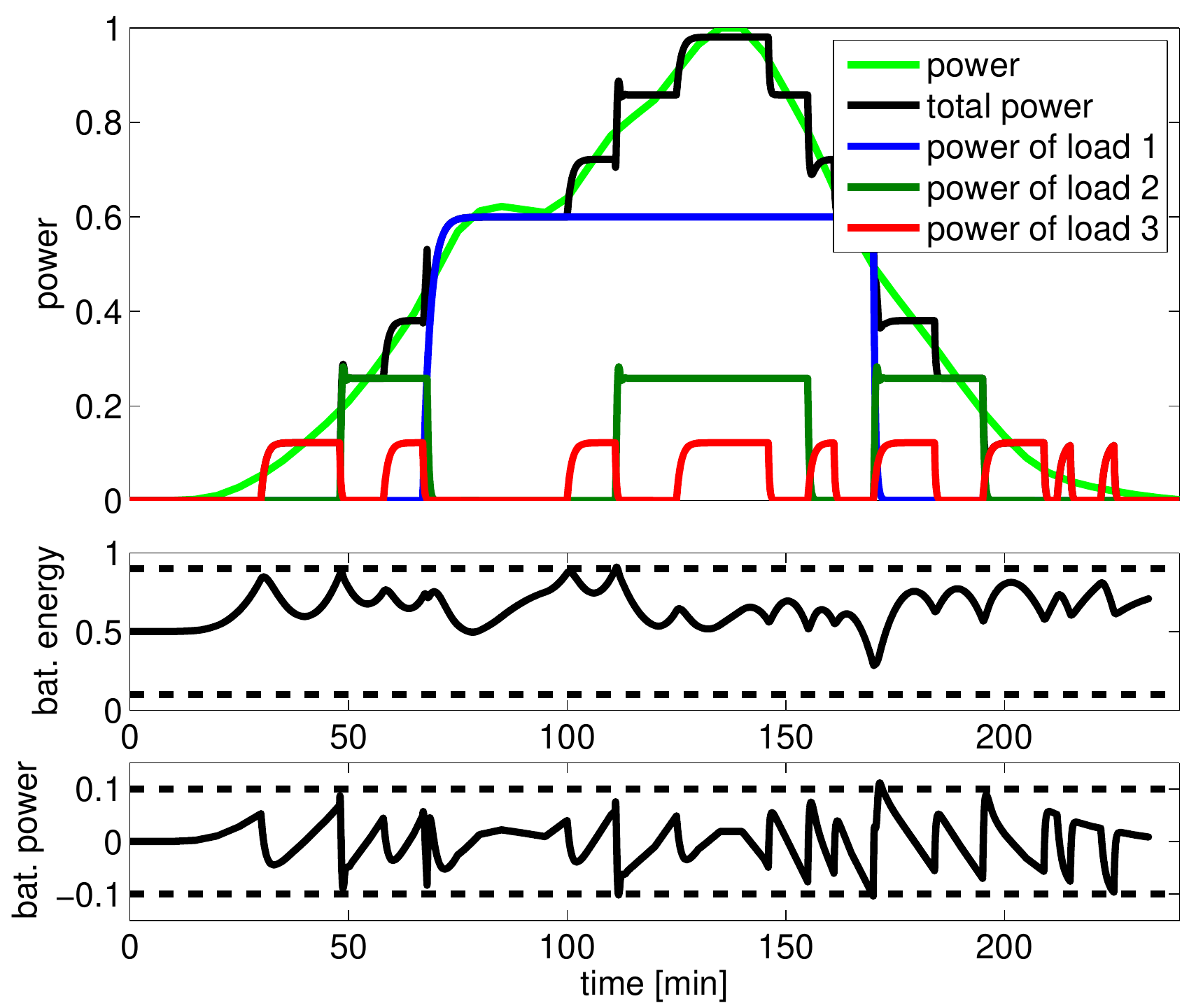}
\caption{Load scheduling with a battery initialized at 50\% SOC. Top figure: Power curve (PU) to be tracked (green line) with total real power demand due to load switching (black line) and individual dynamic load demands (colored lines). Middle and bottom figure are battery energy (100\% SOC) and battery power demand (PU) where a positive value indicates charging.}
\label{bat05}
\end{figure}

For the first results depicted in this paper, the scheduling algorithm is initialized for a battery with an initial SOC of 50\% and leads to the final result depicted in Figure~\ref{bat05}. It can be seen that the load scheduling manages to track the irregular power curve (green line) by scheduling of the loads at the appropriate times, while keeping the battery energy in SOH between the boundaries of 10\% and 90\%. Power demand on the battery also remained within 10\%  of the total power demand in PU.

\begin{figure}[ht]
\centering
\includegraphics[width=.85\columnwidth]{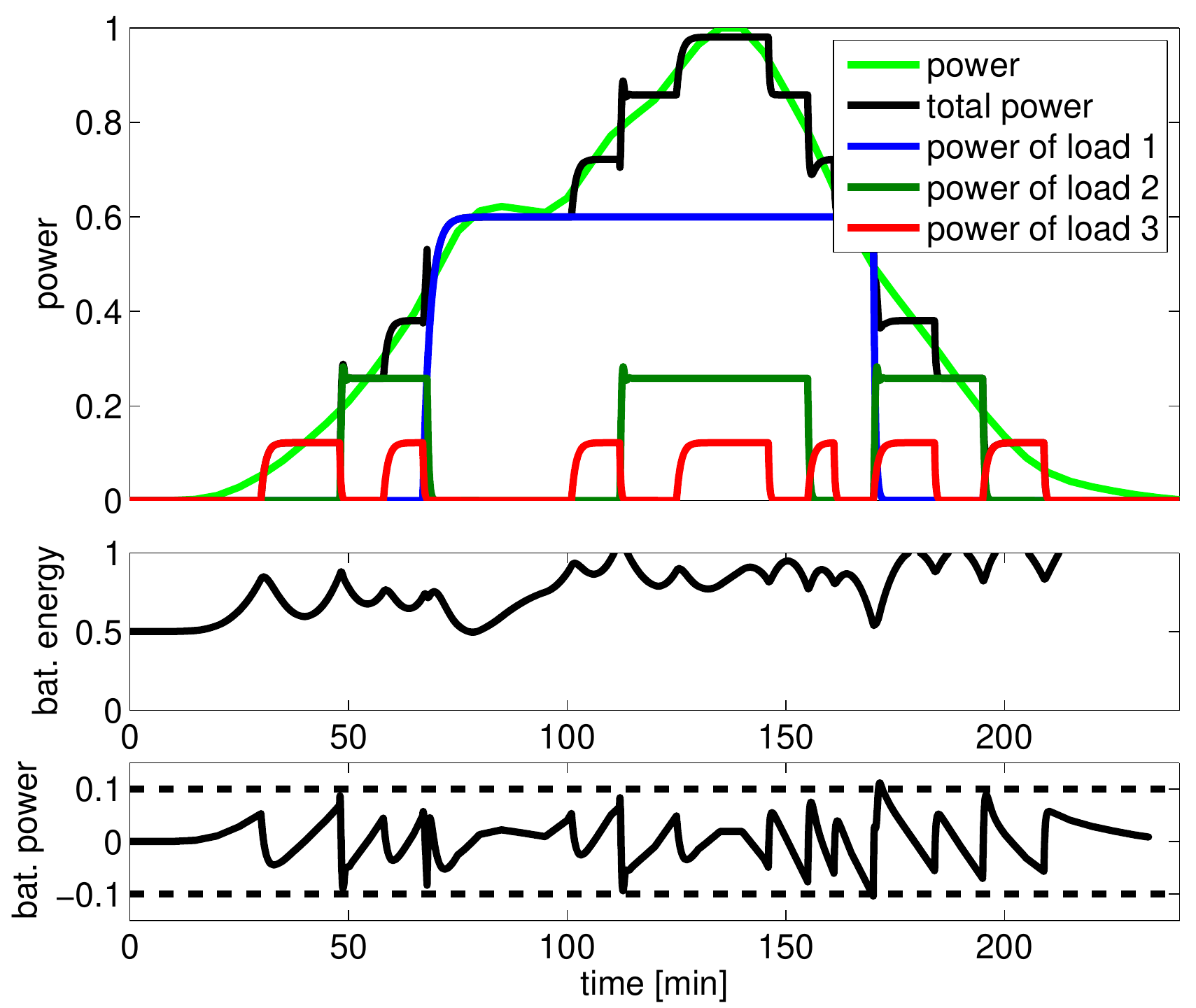}
\caption{Load scheduling with a battery initialized at 50\% SOC, but {\em without\/} constraints on the battery energy. Middle figure, showing the battery energy (100\% SOC), clearly indicates overcharging, despite power curve tracking and constraints on battery power.}
\label{bat05over}
\end{figure}

One may be tempted to concluded that the load switching results are relatively easily obtained due to 50\% SOC level of the battery. However, if the constraints on battery energy weighted by the barrier function $B_j(e(t_m))$ are ignored (e.g. coefficients $C_2 = C_3 =C_4 = 0$), no constraints on battery energy are taken into account. In that case, the battery may be overcharged as indicated in the results summarized in Figure~\ref{bat05over}.

\subsection{Dynamic Load Switching Results with Extreme SOC}

\begin{figure}[ht]
\centering
\includegraphics[width=.85\columnwidth]{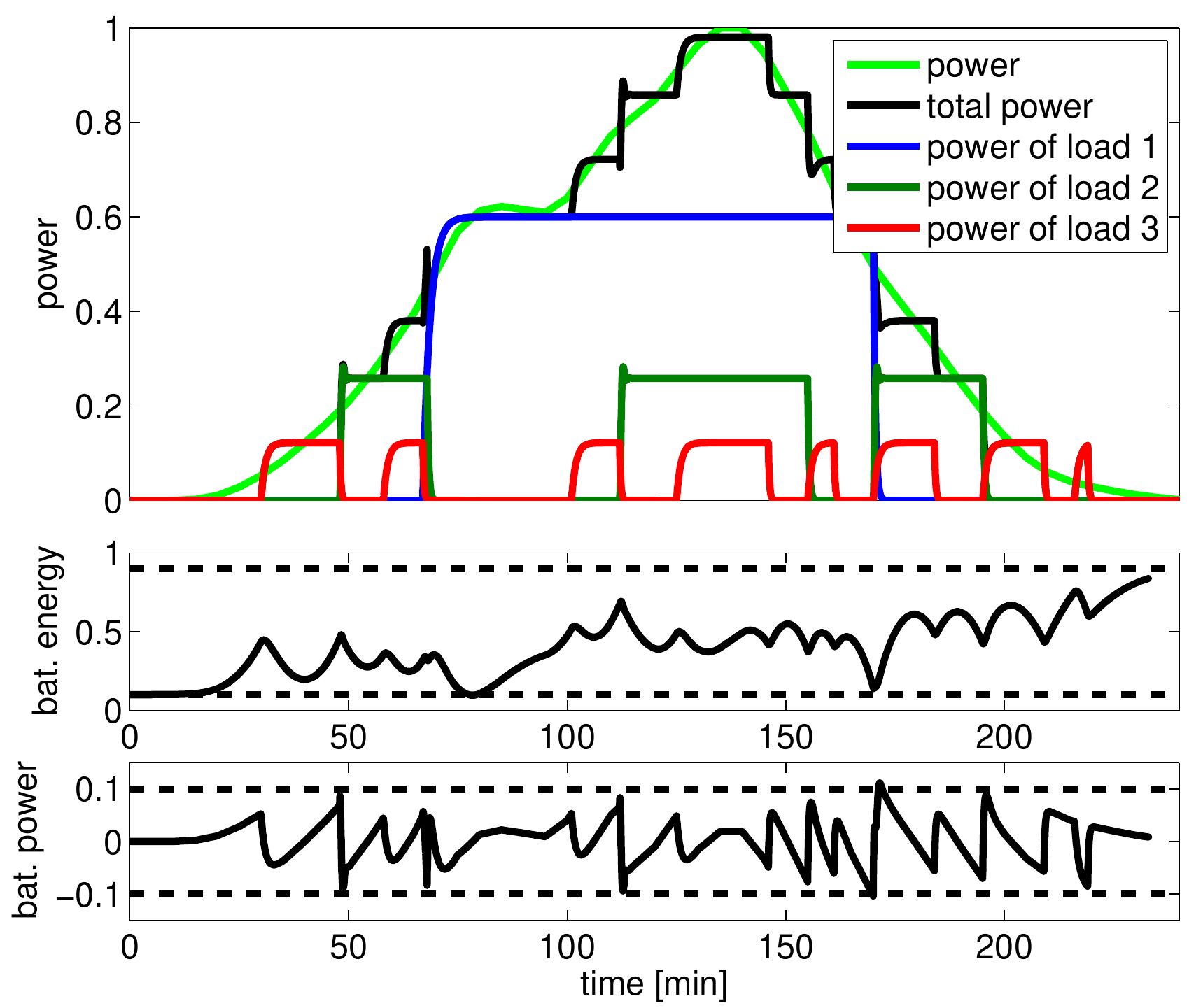}
\caption{Load scheduling with a battery initialized at 10\% SOC. Top figure: Power curve (PU) to be tracked (green line) with total real power demand due to load switching (black line) and individual dynamic load demands (colored lines). Middle and bottom figure are battery energy (100\% SOC) and battery power demand (PU) where a positive value indicates charging.}
\label{bat01}
\end{figure}

\begin{figure}[ht]
\centering
\includegraphics[width=.85\columnwidth]{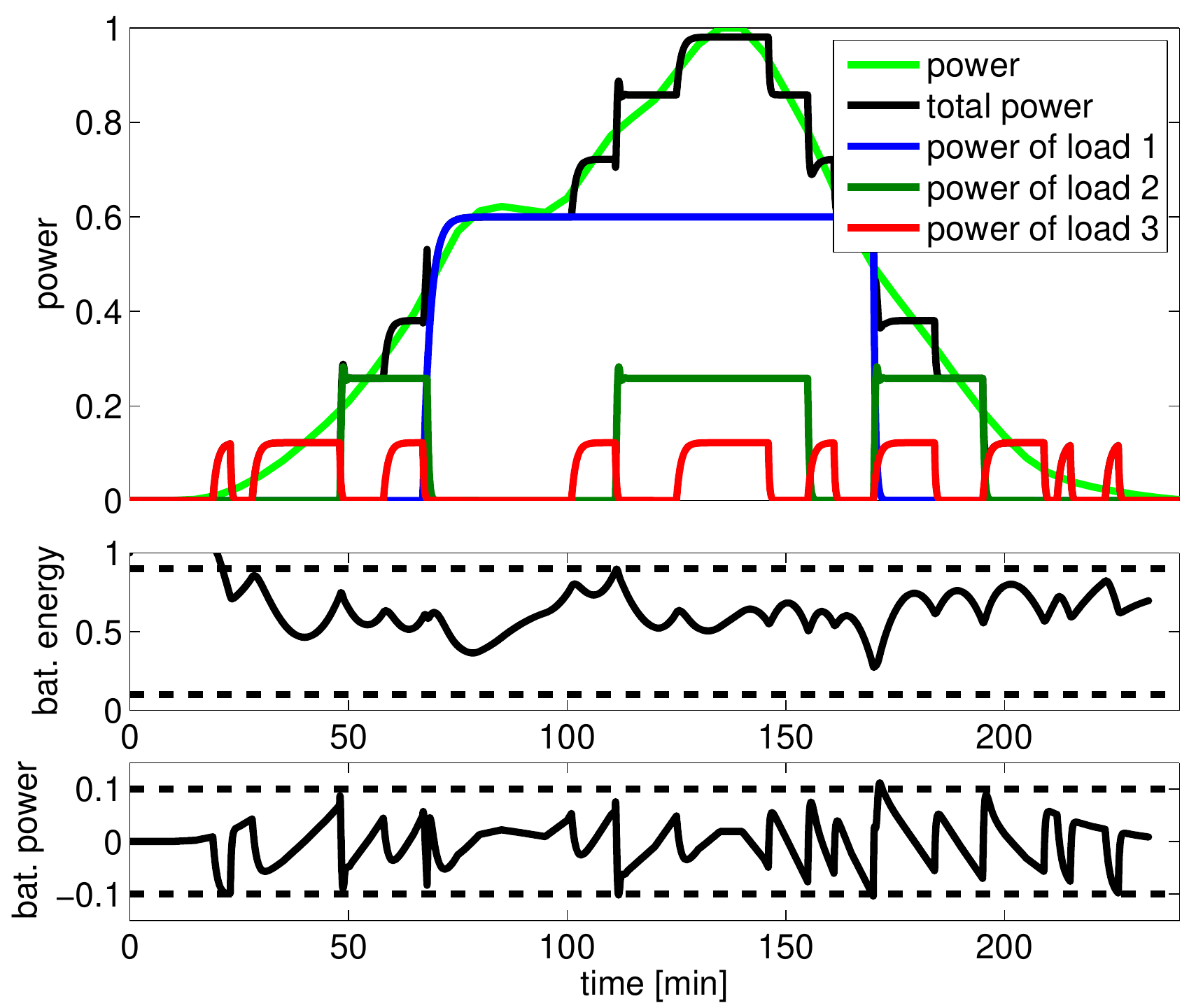}
\caption{Load scheduling with a battery initialized at 100\% SOC. Top figure: Power curve (PU) to be tracked (green line) with total real power demand due to load switching (black line) and individual dynamic load demands (colored lines). Middle and bottom figure are battery energy (100\% SOC) and battery power demand (PU) where a positive value indicates charging.}
\label{bat100}

\end{figure}

Starting the battery at a nearly discharged state (SOC of 10\%) or fully charged state (100\% SOC) would require a careful switching regime of loads at the beginning of the power curve to be tracked. Running the scheduling algorithm for a battery using a moving time prediction horizon where the load scheduling is computed according to the MPC approach \fref{MPC} with the optimization function given in \fref{optcrit} automatically decides on carefully switched loads at the beginning of the power curve to ensure the power curve is tracked, while bringing the battery back into its allowed operating regime. The results of the scheduling algorithm for a battery with a SOC of 10\% and 100\% are summarized in Figure~\ref{bat01} and Figure~\ref{bat100} respectively.

The subtle differences between the load scheduling summarized in Figure~\ref{bat01} and Figure~\ref{bat100} show how the proposed load scheduling algorithm can handle different SOC conditions of the battery, while still tracking the power curve. All this is done, despite the difference in switching dynamics between the loads summarized earlier in Figure~\ref{dynamicsload}.

\section{Conclusions}\label{sec:CCLS}

The optimal scheduling of electrical loads with known and distinct  time dependent power demand profiles is solved by formulating a model predictive approach in which only a short term forecast of power production is needed. The length of the short term forecast is determined by the minimum on/off time of the electrical loads. The optimal scheduling is solved by computing the finite number of possible on/off load switching combinations over the short term power forecast and formulating an objective function that minimizes the difference between (real) power production and (real) power demand, while at the same time taking into account constraints on electrical energy storage and power delivery of a battery system. Simulation results show that relatively short term power forecast profiles can be used to effective schedule dynamic loads with switching dynamics that may be different for each load and can even include load dynamics that have power surge demands.

\bibliographystyle{IEEEtran}
\bibliography{mylib.bib}

\end{document}